\documentclass[11pt]{amsart}

\usepackage{fullpage}

\usepackage[english]{babel}
\usepackage{color}
\usepackage{array,amscd,amsmath,amsthm}
\usepackage{amssymb}
\usepackage[dvips]{graphics}
\usepackage[all]{xy}
\xyoption{2cell}
\allowdisplaybreaks

\newtheorem{theorem}{Theorem}[section]
\newtheorem{proposition}[theorem]{Proposition}

\newtheorem{corollary}[theorem]{Corollary}
\newtheorem{conjecture}[theorem]{Conjecture}

\theoremstyle{definition}
\newtheorem{definition}[theorem]{Definition}
\newtheorem{example}[theorem]{Example}
\newtheorem{remark}[theorem]{Remark}

\newcommand{\Q}{\mathbb{Q}}

\def\M{\mathcal{M}}

\def\Mbar{\overline{\M}}

\def\pa{\partial}

\def\si{\bullet}
\def\Ss{\mathcal{S}}

\def\ra{\rightarrow}
\def\arr#1#2{\stackrel{#1}{#2}}
\def\hra{\hookrightarrow}

\def\A8{A_{\infty}}

\def\l{\lambda}

\def\Uu{\mathcal{U}}

\def\lra{\longrightarrow}

\def\s{\sigma}

\def\ol#1{\overline{#1}}
\def\wh#1{\widehat{#1}}
\def\dis{\displaystyle}

\begin{document}

\title{A remark on the homotopical dimension of
some moduli spaces of stable Riemann surfaces}

\author{Gabriele Mondello}

\address{Massachusetts Institute of Technology - Department of Mathematics
\newline
77 Massachusetts Avenue, Cambridge, MA 02139}

\email{gabriele@math.mit.edu}


\begin{abstract}
Using a result of Harer, we prove certain
upper bounds for the homotopical/cohomological
dimension of
the moduli spaces of Riemann surfaces
of compact type, of Riemann surfaces with rational tails and of Riemann
surfaces with at most $k$ rational components. These bounds would
follow from conjectures of Looijenga and Roth-Vakil.
\end{abstract}


\maketitle

\begin{section}*{Introduction}
Tautological classes, Chow groups and rational cohomology of
the moduli spaces of pointed Riemann surfaces have curious vanishing
properties, which are not completely understood yet.
It seems natural to look for a geometric explanation of these
vanishings, namely these moduli spaces ``should'' be covered by
a certain number of open affine sets, or at least ``should'' have
an affine stratification with a certain number of strata.
Notice that this number also bounds the maximal dimension of a
complete subvariety contained in such a moduli space.

Historically, first Arbarello
\cite{arbarello:weierstrass} realized that a stratification of
the moduli space $\M_g$ of Riemann surfaces of genus $g$
could be useful to investigate the geometric properties
of $\M_g$, starting from the locus of hyperelliptic curves and
then climbing the $g-1$ layers of his stratification.
Ten years later, Diaz \cite{diaz:complete} used a variant of Arbarello's
stratification to show that $\M_g$ cannot contain complete subvarieties
of dimension $g-1$. Some years later, using similar stratifications,
Looijenga \cite{looijenga:tautological} proved that the tautological
classes in the Chow ring of $\M_g$ vanish in degree $\geq g-1$, which implies
Diaz's result. He also showed that the tautological classes
of the moduli space $\M^{rt}_{g,n}$ of $n$-pointed
stable Riemann surfaces of genus $g$ with rational tails
vanish in degree $\geq g-1+n$. Then he formulated the following conjecture.

\begin{conjecture}[Looijenga]\label{conj:looijenga}
Let $g,n\geq 0$ be such that $2g-2+n>0$.\\
The coarse space $M_{g,n}$ of $\M_{g,n}$ has a stratification
\[
M_{g,n}=\coprod_{i=1}^{g-\delta_{n,0}} S_i\qquad\text{with}\quad\ol{S}_j=\coprod_{i\leq j} S_i
\quad\text{for all $j$}
\]
such that each locally closed stratum $S_i$ is affine.
\end{conjecture}

The conjecture would be a consequence of the existence of a
Zariski open cover of $\M_{g,n}$
made of $g-\delta_{n,0}$ affines. On the other hand,
one may notice that both Arbarello's and Diaz's stratifications of $\M_g$ have the
right number of strata, while not even conjectural candidates
have been proposed for the affines.

We stress that Conjecture~\ref{conj:looijenga} has much stronger consequences,
as it is clear from the following general statement.

\begin{proposition}\label{prop:vanish}
Let $M$ be a scheme of finite type of dimension $d$ over the field $\mathbb{K}$ and
let $M=\coprod_{i=1}^k S_i$ be a stratification, with
$\ol{S}_j=\coprod_{i\leq j} S_i$ for all $j$, such that
each locally closed stratum $S_i$ is affine. 
Then $H^r_{\acute{e}t}(M;\mathcal{F})=0$ for all $r>d+k-1$ and all $\mathcal{F}$
constructible (or torsion, or $l$-adic) sheaf on $M$.
Moreover, if $M$ is reduced over $\mathbb{C}$, then the associated topological
space $M^{top}$ is homotopy equivalent
to a CW-complex of dimension at most $d+k-1$.
\end{proposition}

We refer to \cite{roth-vakil:stratification} for a proof of the previous proposition
and for an extensive discussion on affine stratification and covering numbers
and their applications to the moduli spaces of Riemann surfaces.

Looijenga's conjecture can be extended to other partial compactifications
of the moduli space $\M_{g,n}$ of $n$-pointed Riemann surfaces of genus $g$.
The most interesting cases are the moduli space $\M^{ct}_{g,n}$
of stable Riemann surfaces of compact type (i.e. whose dual graph is a tree),
the moduli space $\M^{rt}_{g,n}$ of stable Riemann surfaces with rational
tails, the moduli space $\M^{irr}_{g,n}$ of irreducible stable Riemann
surfaces and the moduli space $\M^{\leq k}_{g,n}$ of stable Riemann surfaces
with at most $k$ rational components.

The interest for the last moduli space $\M^{\leq k}_{g,n}$ is related to a
vanishing result due to Ionel \cite{ionel:topological} in cohomology and to
Graber and Vakil \cite{graber-vakil:star} in the Chow ring: it
says that tautological classes of
$\M^{\leq k}_{g,n}$ vanish in degree
$\geq g+k$, which implies Looijenga's vanishing and Getzler's
conjecture \cite{getzler:virasoro} (which says that $\psi_{i_1}\cdots\psi_{i_g}=0$ on $\M_{g,n}$).
As a consequence of these investigations, Roth and Vakil \cite{roth-vakil:stratification}
wondered about
the meaning of this new stratification of $\Mbar_{g,n}$ by number of rational
components and they raised the following question.

\begin{conjecture}[Roth-Vakil]\label{conj:affine}
Let $g,n>0$ such that $2g-2+n>0$.\\
The coarse space of $\M^{\leq k}_{g,n}$
(resp. $\M^{rt}_{g,n}$, $\M^{ct}_{g,n}$)
has an affine stratification made of $g+k$ (resp. $g+n-1$, $2g-2+n$) strata.
\end{conjecture}

The purpose of the present paper is to analyze the homotopy type
of these moduli spaces of Riemann surfaces, proving as corollary
that their rational cohomological dimension is what we would expect
if Conjecture~\ref{conj:affine} were true,
according to Proposition~\ref{prop:vanish}.
Hence, in a certain sense, our results support Conjecture~\ref{conj:affine}.

Our computations heavily rely on
the case of smooth Riemann
surfaces, which
was investigated by Harer
\cite{harer:virtual} in the following way.

The moduli space $\M_{g,n}$ is the quotient of the
Teichm\"uller space $\mathcal{T}_{g,n}$ by the action of the mapping
class group $\Gamma_{g,n}$. As $\mathcal{T}_{g,n}$ is homeomorphic to a
Euclidean space and $\Gamma_{g,n}$ acts with finite stabilizers,
the orbifold-theoretic quotient $\M_{g,n}$ represents $B\Gamma_{g,n}$.
Hence, the cohomology groups of $\M_{g,n}$ are those of $\Gamma_{g,n}$
essentially by definition; so, instead of dealing with the orbifold $\M_{g,n}$,
Harer could work $\Gamma_{g,n}$-equivariantly on $\mathcal{T}_{g,n}$.
In particular, he found a $\Gamma_{g,n}$-equivariant retraction of $\mathcal{T}_{g,n}$
onto a subcomplex of dimension $4g-4+n$ when $g,n>0$
or $n-3$ if $g=0$, thus bounding 
the homotopical dimension of every unramified cover of $\M_{g,n}$ which is
a manifold, and so
the virtual cohomological dimension of $\Gamma_{g,n}$.
Actually, he proved also that $\Gamma_{g,n}$ is a virtual duality
group and so that the bound is sharp.

However, dealing with partial compactifications of $\M_{g,n}$ which are not
$K(G,1)$'s, we cannot directly reduce to equivariantly work on some partial compactification
of the Teichm\"uller space but we would need to put our hands on the orbifolds themselves.
In fact, it turns out that some objects we are interested in
are simply connected, as shown in Appendix~\ref{app:fundamental}.
Moreover, we would need a {\it virtual} (that is, orbifold) analogue of homotopical
and cohomological dimension of an orbispace.

Instead of treating in full detail these technicalities, we prefer to show how the
main result works for the moduli spaces of Riemann surfaces with level structures
$\M_{g,n}^\l$ whose Deligne-Mumford compactification $\Mbar_{g,n}^\l$
is a smooth variety, sketching the needed modifications to make
the argument work for moduli spaces of Riemann surfaces which are orbifolds.

We incidentally remark that
the cohomology with rational coefficients
of any orbifold coincides with that of its coarse space
and that the virtual homotopical dimension of an orbifold
actually bounds the homotopical dimension
of its coarse space, so that one can also read our result just
at the coarse level, with little loss of information.

Our main result (Theorem~\ref{maintheorem})
concerns any partial compactification of the moduli space of Riemann
surfaces obtained by adding some strata of its Deligne-Mumford compactification.
For the partial compactifications mentioned above, it reduces to the following.

\begin{corollary}\label{maincorollary}
Let $g\geq 0$ and $2g-2+n>0$.
The moduli spaces of Riemann surfaces which are irreducible,
of compact type, with rational tails or with at most $k$ rational
components have the following upper bounds
on their virtual homotopical/cohomological dimension.\\
\begin{center}
\begin{tabular}{c|c|c}
Moduli space of stable & virtual homotopical
& virtual cohomological \\
Riemann surfaces & dimension for $n>0$ & dimension for $n=0$ \\
\hline
$\M^{irr}_{g,n}$ & $4g-3+n$ & $4g-3$ \\
$\M^{ct}_{g,n}$ & $5g-6+2n$ & $5g-6$ \\
$\M^{rt}_{g,n}$ with $g>0$ & $4g-5+2n$ & $4g-5$ \\
$\M^{\leq k}_{g,n}$ & $4g-4+n+k$ & $4g-4+k$  \\ 
\end{tabular}
\end{center}
\end{corollary}

The idea is the following:
given a smooth complex variety $X$, a divisor $D\subset X$ with
normal crossings induces a locally closed stratification
by smooth subvarieties.
One can reconstruct the homotopy type of $X$
from the homotopy type of the strata and of natural
torus bundles over the strata (see Proposition~\ref{prop:homotopy}).
Thus, if one knows the
homotopical (resp. cohomological) dimension of each stratum,
then one can produce an upper bound for the homotopical
(resp. cohomological) dimension of $X$ (see Corollary~\ref{cor:dimension}).

Deligne-Mumford's compactification $\Mbar_{g,n}$ 
has a natural divisor with normal crossings
(in the orbifold sense) $\pa\M_{g,n}$, that is the locus
of singular surfaces. Each stratum $\M_\sigma$ of the induced
stratification parametrizes stable Riemann surfaces of topological
type $\sigma$.

The beauty of Deligne-Mumford's compactification is that all open strata
$\M_\sigma$ are isomorphic to products of smaller moduli spaces of smooth
Riemann surfaces, possibly up to the action of a finite group of symmetries.
Hence, Harer's result allows us to conclude.

We would like to stress that the particolar role played by rational components,
which was noticed by Graber and Vakil, is already present in Harer's result.
In fact, for $n>0$, the virtual homotopical dimension of $\M_{g,n}$
is $4g-4+n$ if $g>0$ but $n-3$ (and not $n-4$ but one more!) if $g=0$.
This shift appears to be responsible, at least at a topological level, for the
homotopical dimension to grow according to the number of rational components.
\begin{subsection}{Acknowledgments}
I would like to thank Enrico Arbarello, John Harer, Jason Starr and Ravi Vakil
for valuable discussions on moduli spaces of Riemann surfaces, and
Mark Behrens and Andr\'e Henriques for useful conversations on Cech-Segal
construction and orbispaces.
I also thank the referee for suggesting to use level structures
and Harvey's bordification to simplify the
argument and for useful comments.
\end{subsection}
\end{section}
%
\begin{section}{Smooth varieties and divisors with normal crossings}\label{sec:dimension}
%
%
%
Let $X$ be a smooth complex variety of dimension $N$ and
$D$ a divisor with
simple\footnote{The simplicity of $D$ is
not necessary but the notation becomes lighter in this case.}
normal crossings, that is $D=\cup_{i\in I} \ol{D}_i$
such that $\ol{D}_i$ is a smooth divisor of $X$ and every two components
$\ol{D}_i$ and $\ol{D}_j$ intersect transversely.
For every $J\subset I$, call $\ol{D}_J:=\cap_{j\in J} \ol{D}_j$
and identify $X$ with $\ol{D_\emptyset}$.

We will denote by $\wh{D}_J$
the manifold with corners obtained performing
a real oriented blow-up $b_J:\wh{D}_J\ra \ol{D}_J$ along
$\dis\bigcup_{i\in I\setminus J}\ol{D}_{J\cup\{i\}}$.
We will also write $\wh{X}$ for $\wh{D}_\emptyset$.
Notice that the natural inclusion
$D_J:=\ol{D}_J\setminus\bigcup_{J\subsetneq L} \ol{D}_L\hra \wh{D}_J$ is a homotopy
equivalence.

For every chain $J_0\subsetneq J_1\subsetneq\dots\subsetneq J_k\subseteq I$, define
$\wh{D}_{J_0,\dots,J_k}:=\wh{D}_{J_0}\times_{\ol{D}_{J_0}}\ol{D}_{J_k}$
and $D_{J_0,\dots,J_k}:=\wh{D}_{J_0}\times_{\ol{D}_{J_0}} D_{J_k}$.
The inclusion $D_{J_0,\dots,J_k}\hra\wh{D}_{J_0,\dots,J_k}$ is a homotopy equivalence
and $\wh{D}_{J_0,\dots,J_k}\ra \wh{D}_{J_k}$ can be identified to
the $(S^1)^d$-bundle
\[
\prod_{i\in J_k\setminus J_0} b_{J_k}^*\left( SN(\ol{D}_i/X)\Big|_{\ol{D}_{J_k}}\right)
\]
where $SN(\ol{D}_i/X)$ is the $S^1$-bundle associated to the normal bundle of $\ol{D}_i$
and $d=|J_k\setminus J_0|$.

Now, define the simplicial topological space (without degeneracies) $X_\si$ as
$X_k:=\coprod_{J_0\subsetneq\dots\subsetneq J_k} \wh{D}_{J_0,\dots,J_k}$, with
the following face maps.
For $l=1,\dots,k-1$, the face map
$f_l:\wh{D}_{J_0,\dots,J_k}\lra \wh{D}_{J_0,\dots,\wh{J_l},\dots,J_k}$
is the identity.
The face map $f_k:\wh{D}_{J_0,\dots,J_k}\lra \wh{D}_{J_0,\dots,J_{k-1}}$
is induced by the inclusion $\ol{D}_{J_k}\hra\ol{D}_{J_{k-1}}$,
whereas $f_0:\wh{D}_{J_0,\dots,J_k}\lra\wh{D}_{J_1,\dots,J_k}$
is induced by the projection
\[
\prod_{i\in J_k\setminus J_0} SN(\ol{D}_i/X)\Big|_{\ol{D}_{J_k}}\lra
\prod_{i\in J_k\setminus J_1} SN(\ol{D}_i/X)\Big|_{\ol{D}_{J_k}}.
\]

\begin{proposition}\label{prop:homotopy}
There is a natural homotopy equivalence between $X$ and the
topological realization
$|X_\si|\simeq
\underset{\quad\longleftarrow}{\mathrm{hocolim}}X_\si$.
\end{proposition}

Let us recall that the topological realization $|X_\si|$ of a simplicial
topological space $X_\si$ is obtained from
$\coprod_{i\geq 0} X_i\times\Delta_i$ by identifying
$(x_i,t_0,\dots,t_{k-1},0,t_{k+1},\dots,t_i)\in X_i\times\Delta_i$
with $(f_k(x_i),t_0,\dots,t_{k-1},t_{k+1},\dots,t_i)\in
X_{i-1}\times\Delta_{i-1}$ for every $i\geq 1$ and $0\leq k\leq i$.

\begin{proof}
Let $\mathcal{I}_d\subset \mathcal{P}(I)$
be the set of subsets $J\subset I$ of cardinality $N-d$ such that
$\ol{D}_J$ is nonempty and let $\mathcal{I}=\bigcup_d\mathcal{I}_d$.

Set $X(-1):=X$ and proceed by induction on $d$.

Let $0\leq d\leq N$ and suppose we have already defined $X(r)$ and
$V_J$ with $J\in\mathcal{I}_r$ for all $0\leq r< d$.
For every $J\in\mathcal{I}_d$ choose a tubular neighbourhood $V_J$ of $\ol{D}_J\cap X(d-1)$ inside
$X(d-1)$ in such a way that the closures of all $V_J$'s inside $X$ are disjoint.
Let $X(d):=X(d-1)\setminus\bigcup_{J\in\mathcal{I}_d} V_J$.

Define the open cover $\mathcal{U}=\{U_J\,|\, J\in\mathcal{I}\}$ of $X$,
with $U_J$ consisting of a small fattening of $\ol{V}_J$ inside $X$, in such a way that:
\begin{itemize}
\item[-]
every $U_J$ retracts by deformation onto $V_J$
\item[-]
$U_J\cap U_K$ is nonempty if and only if $J\subseteq K$ or $K\subseteq J$
\item[-]
for every chain $J_0\subsetneq J_1\subsetneq\dots\subsetneq J_k$
the natural map $U_{J_0}\cap\dots\cap U_{J_k}\arr{\simeq}{\lra} D_{J_0,\dots,J_k}$
is a homotopy equivalence.
\end{itemize}
As the realization of a simplicial topological space represents
the homotopy colimit of the associated system,
the morphism $U_\si \lra X_\si$ defined above
induces a homotopy equivalence $|U_\si|\arr{\simeq}{\lra}|X_\si|$.

By the Cech-Segal construction \cite{segal:classifying},
every partition of unity on $X$ subordinated to the open cover
$\mathcal{U}$ gives a homotopy equivalence
$X \lra |U_\si|$.
\end{proof}

\begin{remark}
If $D$ still has normal crossings but some components of $D$
self-intersect, then the result above holds with minor
(but obvious, and mostly notational) modifications.
\end{remark}

The {\it homotopical dimension} $\mathrm{htdim}(X)$
of a topological space $X$ is
the smallest dimension
of a CW-complex that is homotopically equivalent
to $X$. The {\it cohomological dimension} $\mathrm{chldim}(X)$ of $X$
is the smallest $k\in\mathbb{N}$ such that $H^n(X;\mathbb{L})=0$
for every $n>k$ and every local system $\mathbb{L}$ over $X$.
Clearly, $\mathrm{htdim}(X)\geq\mathrm{chldim}(X)$.

\begin{corollary}\label{cor:dimension}
Let $X$ be a complex manifold and $D$ a divisor with normal crossings.
\begin{itemize}
\item[(a)]
The homotopical dimension $\mathrm{htdim}(X)$ is bounded above
by the maximum
of $\mathrm{htdim}(D_J)+2|J|$ for $J\subseteq I$.
\item[(b)]
The cohomological dimension $\mathrm{chldim}(X)$ is bounded above
the maximum of $\mathrm{chldim}(D_J)+2|J|$
for $J\subseteq I$.
\end{itemize}
\end{corollary}
\begin{proof}
As $|X_\si|$ is a quotient of $\coprod_k \Delta_k\times X_k$,
the space $X$ is homotopy equivalent to a CW-complex made of cells of dimension
at most $\mathrm{max}\{k+\mathrm{htdim}(X_k)\,|\,X_k\neq\emptyset\}$.
Moreover, $\mathrm{htdim}(X_k)\leq \mathrm{htdim}D_{J_0,\dots,J_k}
\leq |J_k|-|J_0|+\mathrm{htdim}D_{J_k}$ and the maximum of
$k+|J_k|-|J_0|+\mathrm{htdim}D_{J_k}$ is attained when $|J_k|=k$ and $|J_0|=0$.

For the cohomological dimension, consider the open cover $\Uu$ as in the proof
of Proposition~\ref{prop:homotopy} (with minor but obvious modifications if
some components of $D$ self-intersect).
Given a local system $\mathbb{L}$ on $X$,
we can build a Leray spectral sequence associated to the
open cover $\Uu$ of $X$ with $E_1$-term
\[
E^{p,q}_1=\bigoplus_{J_0\subsetneq\cdots\subsetneq J_p\in\mathcal{I}}
H^q(U_{J_0}\cap\dots
\cap U_{J_p};\mathbb{L})
\]
converging to $H^{p+q}(X;\mathbb{L})$,
where every component of $d_1$ is induced by the natural inclusion.
The same computation as above gives the wanted upper bound.
\end{proof}

\begin{remark}\label{remark:orbifold}
An orbifold $X$ has {\it virtual homotopical dimension} $\leq d$
if $X$ is homotopy equivalent to an orbisimplicial complex (that is, a complex
of groups) of dimension $d$. Moreover, $X$ has {\it virtual cohomological
dimension} $\leq d$ if $H^k(X;\mathbb{L})=0$ for all $k>d$ and all local
systems $\mathbb{L}$ on $X$.
Statements analogous to Proposition~\ref{prop:homotopy}
and Corollary~\ref{cor:dimension} hold in the orbifold category.
\end{remark}
\end{section}
%
%
\begin{section}{Moduli spaces of Riemann surfaces}
Let $g$ and $n$ be nonnegative integers such that $2g-2+n>0$.
Consider the moduli space
$\M_{g,n}$ of Riemann surfaces
$\Sigma$ of genus $g$ together with $n$ distinct
points $\{p_1,\dots,p_n\}\hra \Sigma$.

Let $\Mbar_{g,n}$ be the Deligne-Mumford compactification of $\M_{g,n}$,
obtained adding Riemann surfaces $\Sigma$ with nodes such that
the marked points
sit in the smooth locus of $\Sigma$ and every rational component of $\Sigma$
contains at least three special points, that is points which are either
marked or nodes.

Both moduli space $\M_{g,n}$ and $\Mbar_{g,n}$ are
orbifolds and $\Mbar_{g,n}$ is compact.
The {\it boundary} $\pa\M_{g,n}:=\Mbar_{g,n}\setminus\M_{g,n}$
is a normal crossing divisor in the orbifold sense.

\begin{definition}
Let $\mathcal{B}_{g,n}$ be the set of all homeomorphism
classes of stable $n$-pointed Riemann surfaces of
(arithmetic) genus $g$.
For every $\sigma\in\mathcal{B}_{g,n}$
the {\it stratum} $\M_\sigma$ of $\Mbar_{g,n}$
is the locally closed suborbifold of Riemann surfaces
that belong to $\sigma$.
\end{definition}

A {\it nonempty open set of strata} of $\Mbar_{g,n}$ is
a nonempty subset $\Ss\subseteq\mathcal{B}_{g,n}$
that satisfies the following condition: if
$\sigma\in\Ss$ is obtained from
$\sigma'\in\mathcal{B}_{g,n}$
pinching some simple closed curves of $\sigma'$ to nodes,
then $\sigma'\in\Ss$.

A nonempty open set of strata $\Ss\subseteq\mathcal{B}_{g,n}$
determines a partial compactification $\Mbar_\Ss:=\bigcup_{\s\in\Ss}\M_\s$ of $\M_{g,n}$
made of strata of $\Mbar_{g,n}$ and this correspondence is bijective.
Clearly, $\mathcal{B}_{g,n}$ corresponds to $\Mbar_{g,n}$ and
the trivial system consisting of the smooth surfaces only corresponds to $\M_{g,n}$.

\begin{example}
The moduli spaces
of irreducible Riemann surfaces $\M^{irr}_{g,n}$,
of Riemann surfaces with rational tails (that is, with an irreducible component
of geometric genus $g$) $\M^{rt}_{g,n}$,
of stable Riemann surfaces with at most $k$ rational components
$\M^{\leq k}_{g,n}$ and
of Riemann surfaces of compact type (that is, whose dual graph is a tree) $\M^{ct}_{g,n}$
are partial compactifications of $\M_{g,n}$ of the type described above.
\end{example}

We would like to apply Corollary~\ref{cor:dimension} to $\Mbar_\Ss$
with $\Ss$ a nonempty open set of strata, but $\Mbar_\Ss$ is usually an orbifold and
not a manifold. To get around this problem we could extend
the treatment of Section~\ref{sec:dimension}
to orbifolds (see Remark~\ref{remark:orbifold}),
which is very natural but requires some technicalities.

Instead, we will treat the case of moduli spaces of
Riemann surfaces with finite level structures $\M_{g,n}^\l$
which are smooth manifolds and
admit a smooth compactification $\Mbar_{g,n}^\l$, and
we will sketch how to deduce the result for $\Mbar_\Ss$.

We denote by $\M_{g,n}^\l\ra \M_{g,n}$ the covering space
associated to a finite index subgroup $\Gamma^\l\subset\Gamma_{g,n}
=\pi_0\mathrm{Diff}_+(\Sigma,\{p_1,\dots,p_n\})$.
It is known since Serre that
$\M_{g,n}^\l$ is a manifold if $\Gamma^\l$ 
acts trivially on $H^1(\Sigma;\mathbb{Z}/\ell\mathbb{Z})$
for some $\ell\geq 3$.
It follows from Looijenga \cite{looijenga:prym}, Pikaart-de Jong 
\cite{pikaart-dejong} and Boggi-Pikaart \cite{boggi-pikaart} that
there exist finite index subgroups $\Gamma^\l\subset\Gamma_{g,n}$ such that
the normalization $\Mbar_{g,n}^\l$ of $\Mbar_{g,n}$ in the field of functions
of $\M_{g,n}^\l$ is a smooth manifold. We refer to these papers and to
\cite{deligne-mumford:irreducibility} for a
detailed treatment of level structures, and to \cite{avc-twisted} for a modular
interpretation of the compactified spaces.

Fix a subgroup $\Gamma^\l\subset\Gamma_{g,n}$ of finite index
such that $\Mbar_{g,n}^\l$ is a smooth manifold.
The boundary $\pa\M_{g,n}^\l:=\Mbar_{g,n}^\l\setminus\M_{g,n}^\l$ is a divisor with
normal crossings and the map $f^\l:\Mbar_{g,n}^\l\ra\Mbar_{g,n}$ is
a finite ramified covering.

A stratum of $\Mbar_{g,n}^\l$ is an irreducible component of
$(f^\l)^{-1}(\M_\sigma)$, where $\sigma\in\mathcal{B}_{g,n}$.
This defines a natural map $\mathcal{B}_{g,n}^\l
\ni \tilde\sigma\mapsto [\tilde\sigma]\in \mathcal{B}_{g,n}$,
where $\mathcal{B}_{g,n}^\l=\{\tilde\sigma\,|\,\M_{\tilde\sigma}
\ \text{is a stratum of}\ \Mbar_{g,n}^\l\}$.
An open set of strata $\Ss^\l\subseteq\mathcal{B}_{g,n}^\l$ is a nonempty subset
such that $\Mbar_{\Ss^\l}:=\bigcup_{\tilde\sigma\in\Ss^\l}\M_{\tilde\sigma}$ is
open inside $\Mbar_{g,n}^\l$.

The first important fact is that
the restriction $\M_{\tilde\sigma}\ra \M_{[\tilde\sigma]}$
is a finite unramified cover.
The second (much more) important fact is Harer's
fundamental result.

\begin{theorem}[\cite{harer:virtual}]
Let $h,m\geq 0$ such that $2h-2+m>0$.
There exists a $\Gamma_{h,m}$-equivariant
deformation retraction of the Teichm\"uller space
$\mathcal{T}_{h,m}$ onto a CW-complex of dimension $4h-4+m$ if $h,m>0$
and of dimension $m-3$ if $h=0$ and $m\geq 3$.
If $m=0$ and $h\geq 2$,
then $\Gamma_h$ has virtual cohomological dimension $4h-5$.
\end{theorem}

We are now ready to prove the following.

\begin{theorem}\label{theorem:level}
Let $g,n\geq 0$ such that $2g-2+n>0$
and let $\Gamma^\l\subset\Gamma_{g,n}$ be a finite index subgroup
such that $\Mbar_{g,n}^\l$ is a smooth manifold.
For every nonempty open set of strata $\Ss^\l\subseteq\mathcal{B}_{g,n}^\l$:
\begin{itemize}
\item[-]
if $n>0$, then the moduli space $\Mbar_{\mathcal{S}^\l}$ has homotopical
dimension at most $4g-4+n+rc(\mathcal{S}^\l)$
\item[-]
if $n=0$ and $\Mbar_{\mathcal{S}^\l}\neq\M_g^\l$,
then the moduli space $\Mbar_{\mathcal{S}^\l}$ has
cohomological dimension at most $4g-4+n+rc(\mathcal{S}^\l)$
\end{itemize}
where $rc(\mathcal{S}^\l)$ is the maximum number of rational components
occurring in $[\tilde\sigma]$ for all $\tilde\sigma\in\mathcal{S}^\l$.
\end{theorem}

\begin{proof}
The homotopical dimension of $\M_{\tilde\sigma}$ can be easily computed
using Harer's result, because $\M_{\tilde\sigma}$ is a finite unramified cover
of $\M_{[\tilde\sigma]}$ and $\M_{[\tilde\sigma]}$ looks like
$\prod_i \M_{h_i,m_i}$
modulo a finite group of symmetries:
for instance, one can proceed inductively on the codimension of the stratum.
The upshot is that, if $\tilde\sigma$ has at least a node or a puncture,
$\mathrm{htdim}(\M_{\tilde\sigma})=4g-4+n-2\nu(\tilde\sigma)+rc(\tilde\sigma)$,
where $\nu(\tilde\sigma)$ is the number of nodes of $[\tilde\sigma]$ and
$rc(\tilde\sigma)$ is the number of rational components of $[\tilde\sigma]$.
The result follows from Corollary~\ref{cor:dimension}.
\end{proof}

\begin{remark}
Notice that $\wh{\M_{g,n}^\l}$ is the quotient $\wh{T}_{g,n}/\Gamma^\l$
of Harvey's bordification $\wh{T}_{g,n}$ of the Teichm\"uller space
(see \cite{harvey:bordification}).
\end{remark}

To deduce a similar result for $\Mbar_\Ss$ (or for any orbifold $\Mbar_{\Ss^\mu}$,
corresponding to a finite level structure $\mu$)
we could invoke Remark~\ref{remark:orbifold}.

At the very end, the whole thing boils down to the following.
Let $\Gamma^\mu\subseteq\Gamma_{g,n}$ be a finite index subgroup and let
$\Ss^\mu\subseteq\mathcal{B}_{g,n}^\mu$ be a nonempty open set of strata.
Pick a finite index subgroup $\Gamma^\l\subset\Gamma^\mu$ such that
$\Mbar_{g,n}^\l$ is a smooth manifold and let $\Ss^\l$ be the set of
strata of $\Mbar_{g,n}^\l$ that map to $\Mbar_{\Ss^\mu}$ through
the natural map $\Mbar_{g,n}^\l\ra\Mbar_{g,n}^\mu$.

Proposition~\ref{prop:homotopy} applied to $X=\Mbar_{\Ss^\l}$
and $D=\Mbar_{\Ss^\l}\cap\pa\M_{g,n}^\l$ tells us that
$\Mbar_{\Ss^\l}$ is homotopy equivalent to $|(\Mbar_{\Ss^\l})_\si|$.
The result for $\Mbar_{\Ss^\mu}$ holds 
(replacing ``dimensions'' by ``virtual dimensions'')
if we can ask this
homotopy equivalence to be $\Gamma^\mu$-equivariant.

This can be done just pulling back the open cover
$\mathcal{U}$ of $\Mbar_{\Ss^\l}$ and the partition
of unity subordinated to $\mathcal{U}$
in the proof of Proposition~\ref{prop:homotopy}
from their analogues on $\Mbar_{\Ss^\mu}$.

Furthermore, one can observe that the singularities of any
finite $\Mbar_{g,n}^\l$
are well-behaved with respect to the stratification induced by $\pa\M_{g,n}^\l$,
so that we can even drop the requirement of $\Mbar_{g,n}^\l$ being
an orbifold.

Hence, we can refine the statement of Theorem~\ref{theorem:level}.

\begin{theorem}\label{maintheorem}
Let $g,n\geq 0$ such that $2g-2+n>0$.
Let $\Gamma^\mu\subset\Gamma_{g,n}$ be a finite index subgroup
and $\Ss^\mu\subseteq\mathcal{B}_{g,n}^\mu$ a nonempty set of strata
such that $\Mbar_{\Ss^\mu}$ is open in $\Mbar_{g,n}^\mu$. Then:
\begin{itemize}
\item[-]
if $n>0$, then the moduli space $\Mbar_{\mathcal{S}^\mu}$ has virtual
homotopical dimension at most $4g-4+n+rc(\mathcal{S}^\mu)$
\item[-]
if $n=0$ and $\Mbar_{\mathcal{S}^\mu}\neq\M_g^\mu$,
then the moduli space $\Mbar_{\mathcal{S}^\mu}$ has
virtual cohomological dimension at most $4g-4+n+rc(\mathcal{S}^\mu)$.
\end{itemize}
\end{theorem}

\begin{corollary}
Let $g\geq 0$ and $2g-2+n>0$.
The moduli spaces of Riemann surfaces which are irreducible,
of compact type, with rational tails or with at most $k$ rational
components have the following upper bounds
on their virtual homotopical/cohomological dimension.\\
\begin{center}
\begin{tabular}{c|c|c}
Moduli space of stable & virtual homotopical
& virtual cohomological \\
Riemann surfaces & dimension for $n>0$ & dimension for $n=0$ \\
\hline
$\M^{irr}_{g,n}$ & $4g-3+n$ & $4g-3$ \\
$\M^{ct}_{g,n}$ & $5g-6+2n$ & $5g-6$ \\
$\M^{rt}_{g,n}$ with $g>0$ & $4g-5+2n$ & $4g-5$ \\
$\M^{\leq k}_{g,n}$ & $4g-4+n+k$ & $4g-4+k$ \\ 
\end{tabular}
\end{center}
\end{corollary}

The same result for Riemann surfaces with rational tails
could be obtained in a simpler way. In fact, it is trivial
for $g=0$ and it reduces
to Harer's result for $n=0,1$.
For $g>0$ and $n>1$, one can notice
that the forgetful map $\M^{rt}_{g,n}\ra\M^{rt}_{g,1}$
is a fiber bundle, whose fiber is a smooth manifold of
(real) dimension $2(n-1)$.

Finally, one can also check that for $\Mbar_{g,n}$
one gets the correct bound.
\end{section}
%
%
\appendix
\begin{section}{About the fundamental
group of some moduli spaces}\label{app:fundamental}
In his fundamental paper \cite{harer:virtual}
on the virtual cohomological
dimension of the mapping class group,
Harer was interested in moduli spaces which
are $K(G,1)$'s. In that case, our definition of
virtual homotopical (resp. cohomological) dimension
reduces to the
homotopical (resp. cohomological) dimension of the classifying space
of a torsion-free subgroup of $G$ of finite index.
However, if our moduli spaces are not $K(G,1)$'s, then
this clearer and simpler definition is no longer available.
In this short appendix we check that most partial compactifications
of $\M_{g,n}$ we are interested in are not $K(G,1)$'s.

\begin{proposition}
The moduli space of Riemann surfaces with rational tails
$\M^{rt}_{g,n}$ is a $K(G,1)$ only for
$(g,n)$ equal to $(0,3)$ or $(g,0),(g,1),(g,2)$ and $g>0$.
\end{proposition}
\begin{proof}
In genus $0$, we have $\M^{rt}_{0,n}=\Mbar_{0,n}$
so that it is always simply-connected and never a $K(G,1)$
(except in the trivial case $n=3$).
For $g>0$, we have that $\M^{rt}_{g,n}=\M_{g,n}$ for $n=0,1$
and also that $\M^{rt}_{g,2}$ is the universal curve over
$\M_{g,1}$. For $n>2$ there is a topological fibration
$F\ra\M^{rt}_{g,n}\ra\M_{g,1}$.
Consider $[S,p]\in\M_{g,1}$ and let $F_{[S,p]}$ be the fiber
over $[S,p]$. Then $F_{[S,p]}$ is a projective variety of
complex dimension $n-1$. If $n>2$, then $F_{[S,p]}$ contains
a rational curve, so that $\pi_2(F)\neq 0$ and
$\pi_2(\M^{rt}_{g,n})\neq 0$.
\end{proof}
\begin{proposition}
The moduli space of irreducible stable Riemann surfaces
$\M^{irr}_{g,n}$ is a $K(G,1)$ only if $g=0$ (when
it coincides with $\M_{0,n}$) and it is simply-connected
otherwise.
\end{proposition}
\begin{proof}
Trivially $\M^{irr}_{0,n}=\M_{0,n}$. Instead, for $g\geq 1$
the mapping class group $\Gamma_{g,n}$ is generated by
Dehn twists around nondisconnecting simple closed curves,
so that $\M^{irr}_{g,n}$ is simply connected. However
$H^2(\M^{irr}_{g,n};\Q)$ is not zero for $g>0$
(see~\cite{arbarello-cornalba:algebraic}).
\end{proof}
\begin{proposition}
The moduli space $\M^{\leq k}_{g,n}$ of stable Riemann
surfaces with at most $k$ rational components is a $K(G,1)$
only for $(g,n,k)$ equal to $(0,n,1)$ (when it coincides with $\M_{0,n}$)
and $(1,n,0)$ (when it coincides with $\M_{1,n}$) and it is
simply-connected otherwise.
\end{proposition}
\begin{proof}
The genus $0$ case is trivial because $\M^{\leq 0}_{0,n}$
is empty, $\M^{\leq 1}_{0,n}$ is equal to $\M_{0,n}$ and
$\M^{\leq k}_{0,n}$ is simply connected for $k\geq 2$,
but $H^2(\M^{\leq k}_{0,n};\Q)\neq 0$ for $k\geq 2$.
For $g=1$ and $k=0$ we have $\M^{\leq 0}_{1,n}=\M_{1,n}$.
Instead, for $g=1$ and $k\geq 1$ or for $g>1$ and $k\geq 0$
the moduli space $\M^{\leq k}_{g,n}$ is simply connected
because $\Gamma_{g,n}$ is generated by Dehn twists around
nondisconnecting curves. One can check that in all these
cases $H^2(\M^{\leq k}_{g,n};\Q)$ does not vanish.
\end{proof}
\begin{proposition}
The moduli space $\M^{ct}_{g,n}$ of stable Riemann surfaces
of compact type is a $K(G,1)$ only for $(g,n)$ equal to
$(0,3),(1,1),(1,2),(2,0)$.
\end{proposition}
\begin{proof}
$\M^{ct}_{0,n}$ is $\Mbar_{0,n}$ and $\M^{ct}_{1,n}$
is exactly $\M^{rt}_{1,n}$. An exceptional case is
$\M^{ct}_{2,0}$ which coincides with $\mathcal{A}_2$,
the moduli space of principally polarized Abelian varieties
of dimension $2$. For $(g,n)$ greater or equal to $(2,2)$,
$(3,1)$ or $(4,0)$ the moduli spaces $\M^{ct}_{g,n}$ contain
a complete rational curve, so that $\pi_2(\M^{ct}_{g,n})\neq 0$.
For $(g,n)=(2,1),(3,0)$, one can notice that the maps into
Deligne's torsors (see~\cite{hain-looijenga:moduli})
$\M^{ct}_{2,1}\ra\mathcal{D}_2$ and
$\M^{ct}_{3,0}\ra\mathcal{D}^{pr}_3$
induce isomorphisms on fundamental groups, but
not on rational cohomology. Moreover $\mathcal{D}_g$ and
$\mathcal{D}^{pr}_g$ are $K(G,1)$'s. So $\M^{ct}_{g,n}$
is a $K(G,1)$ only for $(g,n)$ equal to $(0,3)$,
$(1,1)$, $(1,2)$ and $(2,0)$.
\end{proof}
\end{section}
%
%
\bibliographystyle{amsalpha}
\bibliography{virtual-final}
\end{document}